\documentclass[sn-mathphys,Numbered]{sn-jnl}

\usepackage{graphicx}
\usepackage[utf8]{inputenc}
\usepackage{amssymb,amsmath,amsfonts}
\usepackage{mathtools}
\usepackage{xcolor}

\usepackage{multirow}%
\usepackage{amsthm}%
\usepackage{mathrsfs}%
\usepackage[title]{appendix}%
\usepackage{textcomp}%
\usepackage{manyfoot}%
\usepackage{booktabs}%
\usepackage{algorithm}%
\usepackage{algorithmicx}%
\usepackage{algpseudocode}%
\usepackage{listings}%

\newcommand{\ud}{\mbox{d}}
\newcommand{\ba}{\mathbf{a}}
\newcommand{\bb}{\mathbf{b}}
\newcommand{\baf}{\mathbf{a}_f}
\newcommand{\bbf}{\mathbf{b}_f}
\newcommand{\gradphi}{\nabla \phi}
\newcommand{\xf}{\mathbf{x}_f}
\newcommand{\xc}{\mathbf{x}_i}
\newcommand{\xj}{\mathbf{x}_j}
\newcommand{\nf}{\mathbf{n}_f}
\newcommand{\rf}{\mathbf{r}_{f}}
\newcommand{\xfone}{\mathbf{x}_{f,1}}
\newcommand{\xftwo}{\mathbf{x}_{f,2}}
\newcommand{\delfone}{\delta_{f,1}}
\newcommand{\delftwo}{\delta_{f,2}}

\begin{document}

\title[A generalized formulation for gradient schemes in unstructured finite volume method]{A generalized formulation for gradient schemes in unstructured finite volume method}

\author[1]{\fnm{Mandeep} \sur{Deka}}\email{mandeepdeka@iisc.ac.in}

\author[2]{\fnm{Ashwani} \sur{Assam}}\email{aashwani@iitp.ac.in}

\author*[3]{\fnm{Ganesh} \sur{Natarajan}}\email{n.ganesh@iitpkd.ac.in}

\affil[1]{\orgdiv{Chemical Engineering Department}, \orgname{Indian Institute of Science Bangalore}, \orgaddress{\city{Bengaluru}, \postcode{560012}, \state{Karnataka}, \country{India}}}

\affil[2]{\orgdiv{Mechanical Engineering Department}, \orgname{Indian Institute of Technology Patna}, \orgaddress{\city{Patna}, \postcode{801106}, \state{Bihar}, \country{India}}}

\affil*[3]{\orgdiv{Mechanical Engineering Department}, \orgname{Indian Institute of Technology Palakkad}, \orgaddress{\city{Palakkad}, \postcode{678623}, \state{Kerela}, \country{India}}}

\abstract{We present a generic framework for gradient reconstruction schemes on unstructured meshes using the notion of a dyadic sum-vector product. The proposed formulation reconstructs centroidal gradients of a scalar from its directional derivatives along specific directions in a suitably defined neighbourhood. We show that existing gradient reconstruction schemes can be encompassed within this framework by a suitable choice of the geometric vectors that define the dyadic sum tensor. The proposed framework also allows us to re-interpret certain hybrid schemes, which might not be derivable through traditional routes. Additionally, a generalization of flexible gradient schemes is proposed that can be employed to enhance the robustness of consistent gradient schemes without compromising on the accuracy of the computed gradients.}

\keywords{Gradient reconstruction, Generalized formulation, Finite Volume method, Unstructured grid}

\pacs[MSC Classification]{65N08, 76M12}

\maketitle

\section{Introduction}
\label{sec1}

The computation of gradients of flow properties plays a pivotal role in obtaining accurate numerical solutions of Navier-Stokes equations on unstructured meshes. This is because these gradients are necessary to evaluate both inviscid and viscous fluxes in finite volume flow solvers. The solution gradients are employed for higher-order convective schemes for inviscid flux calculations while viscous fluxes are directly dependent on the gradients of velocities and temperature. The construction of robust and accurate flow solvers for incompressible and compressible flows therefore necessitate an accurate and robust methodology for gradient computation among other requirements.

Cell-centered finite volume codes which typically employ unstructured grids reconstruct the centroidal gradients using reconstruction schemes that can be generally classified as belonging to either least squares (LSQ) or Green-Gauss (GG) family of approaches \cite{blazek_computational_2005,moukalled2016finite}. While the former category of techniques are based on the principle of least squares minimization, the latter class of methods rely on the Gauss divergence theorem to calculate the gradients. These reconstruction strategies have been in widespread use, with LSQ approaches largely favoured for Euler flows while GG-based schemes are more common in viscous computations. The schemes of the two families have their own advantages and drawbacks. While the linear least squares approach gives consistent gradients and accurate estimates on isotropic irregular meshes, it is quite inaccurate on stretched meshes over curved geometries \cite{mavriplis2003revisiting}. The use of weights mitigates the problem but the weighted least squares reconstruction requires limiting to stabilise the numerical solution on these meshes even for subsonic inviscid flows \cite{shima2013green} and also results in incorrect gradients on triangulated meshes with skewed elements \cite{mavriplis2003revisiting}. The choice of the explicit neighbourhood necessary for least squares family of schemes affects both the stability and accuracy and requires special attention. 
The Green-Gauss family of schemes do not always lead to consistent gradients but can produce accurate estimates of gradients on curved meshes with high aspect ratio \cite{mavriplis2003revisiting}. Coupled with their ease of implementation, they are a popular choice in viscous flow solvers where they estimate aerodynamic coefficients quite accurately \cite{shima2013green}. The detailed investigations of Diskin and Thomas \cite{diskin2008accuracy} for high Reynolds number applications also demonstrate that the gradient accuracy on high aspect ratio meshes depends both on the grid and the solution. 

It is evident that accuracy, robustness and computational cost are conflicting requirements for any gradient scheme. While there exists no known ``all-in-one'' scheme with all these attributes for a wide range of flow conditions on arbitrary polyhedral meshes, there have been advancements in devising new improved schemes as well as constructing hybrid reconstruction approaches. The latter class of methods attempt to blend the GG and LSQ strategies with a view to jointly harness their advantages and the GLSQ scheme of Shima and co-workers \cite{shima2013green} is a notable effort in this direction. Among the recently proposed schemes for gradient computations are the variational reconstruction (VR) approach of Wang et al. \cite{wang2017compact}, the Implicit Green-Gauss (IGG) \cite{nishikawa2018hyperbolic} scheme, the Taylor-Gauss gradients (TG) \cite{oxtoby2019family,SYRAKOS2023108} and the Modified Green-Gauss (MGG) \cite{deka2018new} reconstruction. The first two schemes are inherently implicit and may be regarded effectively as an extension of compact schemes for unstructured meshes, with the VR scheme categorised as a least squares based method and IGG reconstruction belonging to Green-Gauss family. The TG gradients are constructed similar to least squares but have dependence on the face normal vector similar to GG, thereby having a unique flavor of hybridization. The MGG reconstruction, derived from a different variant of the Gauss-divergence theorem, reconstructs centroidal gradients from face normal derivatives of the scalar which makes it an implicit scheme on non-orthogonal meshes. All the schemes however result in consistent gradients on arbitrary mesh topologies and exhibit advantages compared to their conventional counterparts for fluid flow problems. 
A good overview of various gradient techniques in the finite volume framework including recent developments may be found in \cite{feng2020cell}.
   
The least squares and Green-Gauss reconstruction schemes have been traditionally considered as two different families for gradient computation. This is arguably because they adopt vastly different mathematical philosophies to compute the same quantity of interest. However, efforts to hybridize the constructional philosophies of both classes have yielded some interesting insights into the methods. For instance, it is possible to derive the Green-Gauss formula from a face-based least squares reconstruction using a specific choice of weights \cite{deka2023least}. In another work by Syrakos et al. \cite{SYRAKOS2023108}, a unified framework for gradient reconstruction schemes is proposed which encompasses both GG and LSQ classes of schemes. The basic idea is to use Taylor expansion to express the scalar at a neighbour point (which may or may not be a neighbour cell center) in terms of the value and gradients of the scalar at the cell center where the gradients are sought. The difference in values is then multiplied by a vector and the summation across a suitable stencil yields a system of equations that can be solved to obtain the centroidal gradients. The choice of the neighbour point and the vector used to left multiply the system, distinguishes the different gradients schemes, with LSQ class of methods generated by an orthogonal projection and the GG class of methods by an oblique projection. 

In this work, we propose a generalized formulation for nominally first order accurate gradient schemes in unstructured finite volume method, in the similar spirit as Syrakos et al.\ \cite{SYRAKOS2023108}. Our proposed framework reconstructs centroidal gradients from the directional derivatives of the scalar in a specified neighbourhood, which brings in more generality in the formulation. The mathematical framework of the proposed formulation is presented in section 2. 
Following that in section 3, we demonstrate how existing gradient schemes of both GG and LSQ families can be constructed using the generalized formulation. The framework also enables us to interpret and possible derive hybrid schemes that could retain some benefits of both GG and LSQ family of schemes. This is discussed in section 4. Finally, we extend the proposed framework to construct a generalized version of flexible gradient scheme \cite{nishikawa2021flexible} in section 5. The flexible scheme has a free-parameter that can be tuned to increase robustness of finite volume computations without compromising the gradient reconstruction accuracy.
Some general discussions on the novelties and drawbacks of the proposed formalism are presented in section 6.

\begin{figure}
   \centering
    \includegraphics[scale=0.25]{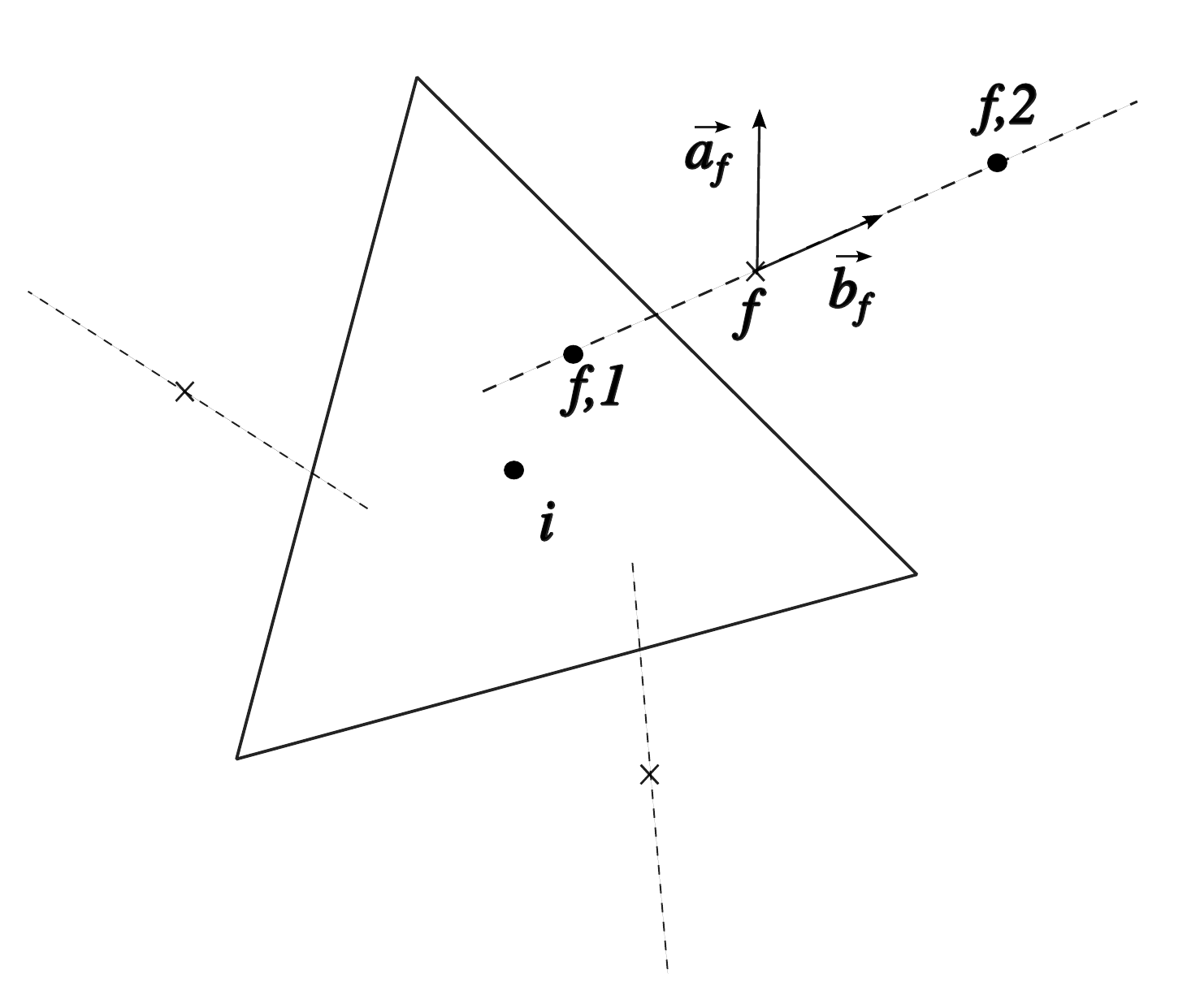}
    \caption{Configuration of a cell in an unstructured mesh with a compact ``neighbourhood'' denoted by a ``$\times$'' symbol. A representative point in the neighbourhood is denoted by ``$f$'' where the two geometric vectors are $\baf$ and $\bbf$ are shown. The points ``$f,1$'' and ``$f,2$'' are along the vector $\bbf$.}
    \label{fig:unified_geo}
\end{figure}

\section{Dyadic sum-vector product formalism}
\label{sec2}

We present a generic mathematical framework for nominally first order accurate gradient reconstruction schemes in unstructured cell-centered finite volume method. {Traditional approaches of deriving gradient schemes are generally through some variants of the Gauss-divergence theorem or through a least squares minimization problem. We, however, adopt a different approach wherein the development of a gradient reconstruction formula is through a tensor identity, that requires the definition of two key entities - the `neighbourhood' of a cell, and, a geometric vector pair ($\ba, \bb$) defined at each element forming the neighbourhood of a cell.} 

To illustrate, let us consider a typical unstructured mesh with a characteristic dimension `$h$'\footnote{For a general unstructured mesh in a $d$-dimensional domain, $h = (V/N)^{1/d}$, where $V$ is the total volume of the computational domain and $N$ is the number of cells discretizing the domain.}, where the gradients of a scalar are required to be computed in a representative cell, with the cell center ``$i$'' (see Fig.\ \ref{fig:unified_geo}).
We define a `neighbourhood' of the cell comprising of a set of points that maintain a compact support about the cell centroid. For the sake of generality, the location of these points are kept arbitrary but typical choices could be (but not restricted to) the cell centres of face-sharing or node-sharing cells, face centers of the faces forming the cell, mid-point of the line joining the face-sharing cell centers, etc. We use the subscript ``$f$'' to denote a representative point in the defined neighbourhood of the cell ``$i$'', as shown in Fig.\ \ref{fig:unified_geo}. With the neighbourhood set, we define a pair of ``geometric'' vectors $\baf$ and $\bbf$ at each point forming the cell neighbourhood. The term ``geometric'' here implies that the vectors shall be function of the mesh metrics alone with a constraint on the magnitudes of these vectors,
\begin{equation}
    |\:\baf\:| \sim O(h^p) \: , \: \: \: \: \: |\:\bbf\:| \sim O(h^q),
\label{eq0}
\end{equation} 
where, $p$ and $q$ are positive integers\footnote{It is presumed that $h \ll 1$ as a requirement for sufficient resolution of the calculations.}. 

For a continuous vector field $\bf{u}$, we can write the following tensor identity,
\begin{equation}
    ({\bf a}_{f} \otimes {\bf b}_{f}) \cdot {\bf u}_f = {\bf a}_f~ ({\bf b}_f \cdot {\bf u}_f) ,
    \label{eq1}
\end{equation}
where ${\bf u}_f$ denotes the value of the vector at the point ``$f$''. Since $\bf u$ is a continuous vector field, we use a Taylor expansion to express ${\bf u}_f$ on the left hand side of Eq.\ \ref{eq1} in terms of the value at the cell centroid ${\bf u}_i$, which after rearranging leads to the expression,
\begin{equation}
     ({\bf a}_f \otimes {\bf b}_f) \cdot {\bf u}_i = {\bf a}_f~ ({\bf b}_f \cdot {\bf u}_f) - ({\bf a}_f \otimes {\bf b}_f) \cdot ({\nabla {\bf u}}_i \cdot ({\bf x}_f - {\bf x}_i) + \epsilon_1) ,
    \label{eq2}
\end{equation}
where, $\epsilon_1 \sim O(h^2)$ represents the higher order terms. Summing up the above equation over the entire neighbourhood of the cell, we get the relation, 
\begin{equation}
   \sum_f ({\bf a}_f \otimes {\bf b}_f) \cdot {\bf u}_i = \sum_f {\bf a}_f~ ({\bf b}_f \cdot {\bf u}_f) -  \sum_f ({\bf a}_f \otimes {\bf b}_f) \cdot ({\nabla {\bf u}}_i \cdot ({\bf x}_f - {\bf x}_i) + \epsilon_1) .
    \label{eq3}
\end{equation}
The sum of the dyadic product of the vectors $\mathbf{a}_f$ and $\mathbf{b}_f$ appearing in the left hand side of the above equation, is a second order tensor, referred to from here on as the `dyadic sum',
 \begin{equation*}
        {\mathbb P} = \sum_f \left({{\bf a}_f} \otimes {{\bf b}_f}\right).
        \label{eq4}
\end{equation*}
To construct the gradient reconstruction scheme, we place another constraint on the choice of the neighbourhood and the geometric vectors, which is that the dyadic sum, $\mathbb{P}$, is invertible. Now, we simply substitute the vector field as the gradient of a scalar field, i.e. $\bf u = \nabla \phi$ in Eq.\ \ref{eq3} and obtain the following expression for centroidal gradients,
\begin{equation}
 \nabla \phi_i = \mathbb{P}^{-1} \sum_f \baf (\bbf \cdot \gradphi_f).
 \label{eq5}
\end{equation}
In obtaining the above expression, the contribution from the second term on the right hand side of Eq.\ \ref{eq3} is ignored. It can be shown that the contribution from that term is $O(h)$ (see appendix \ref{secapp:error_term_1}) and hence can be ignored to obtain an expression for nominally first order accurate gradients. 

The generalized gradient scheme of Eq.\ \ref{eq5} reconstructs centroidal gradients from the directional derivatives of $\phi$ along the $\bbf$ vectors. It is possible to generalize existing gradient reconstruction schemes using this framework, as shall be demonstrated in the next few sections. 
With the neighbourhood and geometric vectors defined, the next crucial step in the formalism is the evaluation of quantity $\bbf \cdot \gradphi_f$ on the right hand side of Eq.\ \ref{eq5}, from which the centroidal gradients are to be reconstructed. 
For a cell-centered finite volume discretization, we have the values of $\phi$ only at the cell centers. Therefore, we can express the quantity $(\bbf \cdot {\gradphi}_f)$, in terms of centroidal values of $\phi$. To do that, we write a general expression for the directional derivative using Taylor series expansion,
\begin{equation}
\begin{aligned}
   {\bf b}_f \cdot \nabla {\phi}_f = \dfrac{\phi(\xftwo) - \phi(\xfone)}{\delftwo - \delfone} - (\delftwo + \delfone) ( {\nabla \nabla \phi}_f : (\bbf \otimes \bbf ) + \epsilon_2  
   \end{aligned}
   \label{eq7}
\end{equation}
where, the $\xfone = \xf + \delfone \bbf$ and $\xftwo = \xf + \delftwo \bbf$ are points along $\bbf$, at a distance of $\delfone/|\bbf|$ and $\delftwo/|\bbf|$ from $\xf$, respectively (see Fig.\ \ref{fig:unified_geo}). The first term on the right hand side of Eq.\ \ref{eq7} essentially expresses $(\bbf \cdot {\gradphi}_f)$ in terms of $\phi$ at specific locations along the $\bbf$ vector. While the choice of the locations $\xfone$ and $\xftwo$ are generic, to ensure compactness, we impose a restriction that $|\delfone\bbf|\sim O(h)$ and $|\delftwo\bbf| \sim O(h)$, i.e., the points $\xfone$ and $\xftwo$ are at an $O(h)$ distance from $\xf$.   
Therefore, substituting Eq.\ \ref{eq7} into Eq.\ \ref{eq5}, we obtain a nominally first order accurate gradient reconstruction formula,
\begin{equation}
\begin{aligned}
 \nabla \phi_i = \mathbb{P}^{-1} \sum_f \baf \left( \dfrac{\phi_{f,2} - \phi_{f,1}}{\delftwo - \delfone} \right),
\end{aligned}
 \label{eq9}
\end{equation}
where, $\phi_{f,1} = \phi(\xfone)$ and $\phi_{f,2} = \phi(\xftwo)$. The contribution from the second term on the right hand side of Eq.\ \ref{eq7} is ignored as it can be shown to be $O(h)$ in magnitude (see appendix \ref{secapp:error_term_2}). 

The expression in Eq.\ \ref{eq9} represents an explicit gradient reconstruction formula if the values of $\phi$ are available at $\xfone$ and $\xftwo$. In general, $\bbf$ need not pass through cell centers, implying that there need not exist a point along the vector $\bbf$ where the values of $\phi$ are available. Therefore, in order to obtain these values, one may need to interpolate. If that is the case, then it can be shown (see appendix \ref{secapp:error_term_3}) that the interpolation needs to be at least second order accurate in order to obtain first order accurate centroidal gradients.
Instead of obtaining $\phi_{f,1}$ and $\phi_{f,2}$, using the values from neighbour cells, one may also use the values and the gradients from neighbour cells in the interpolation scheme, thereby creating an implicit gradient scheme. Implicit schemes can also be constructed by directly expressing the directional derivative $(\bbf \cdot \gradphi_f)$ in terms of the values and gradients in the neighbourhood of the cell. 

\begin{figure}
   \centering
    \includegraphics[scale=0.3]{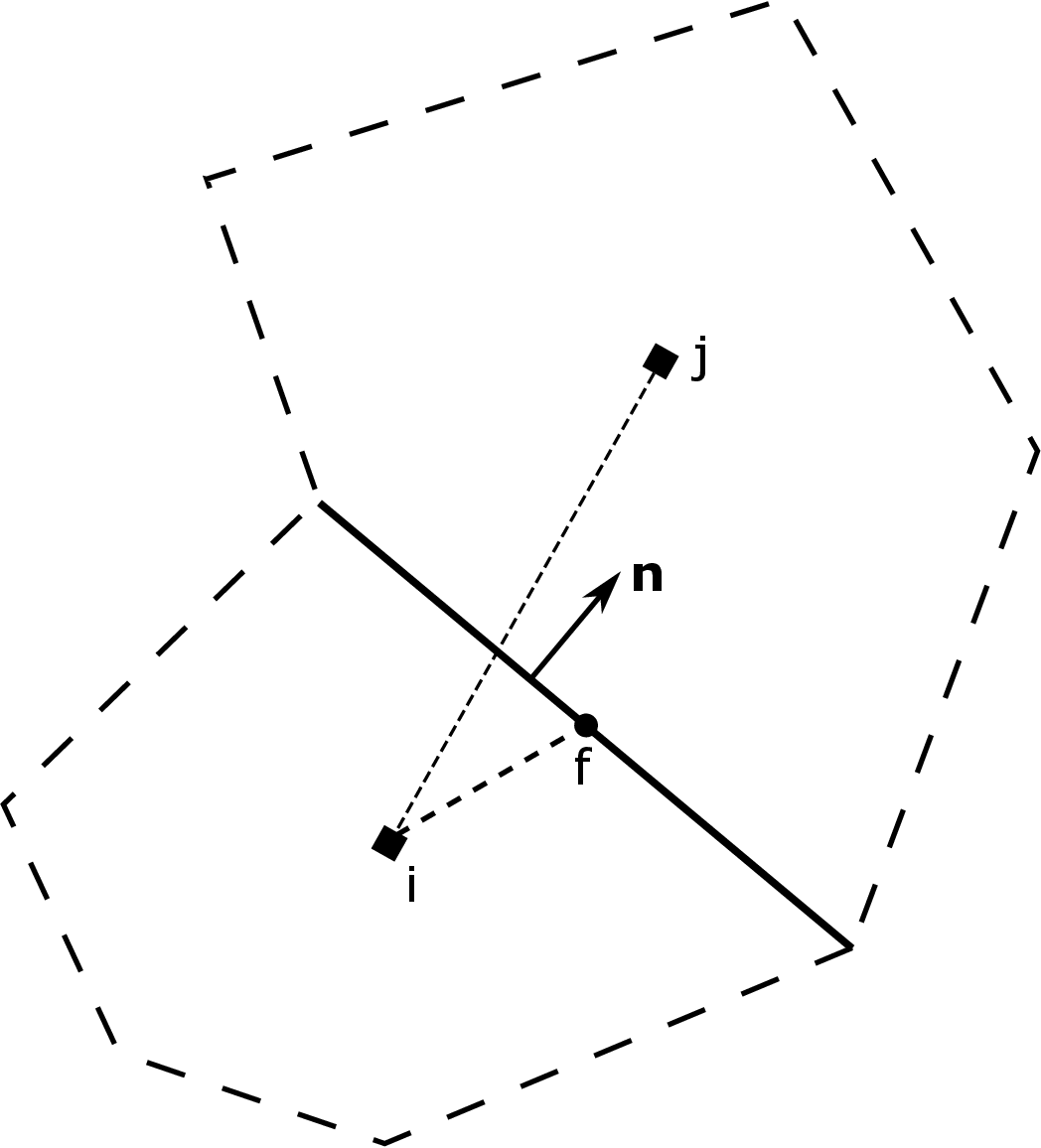}
    \caption{Typical polygonal geometry showing the cells sharing a face. The cells are denoted by their centroids $i$ and $j$, while $f$ represents the face as well as the face center.}
    \label{geo_cv}
\end{figure}

\section{Revisiting the existing gradient reconstruction schemes}

In this section, we demonstrate that existing gradient schemes including the popularly employed Green-Gauss reconstruction \cite{sozer2014gradient,syrakos2017critical} and the linear least squares reconstruction \cite{mavriplis2003revisiting,diskin2011comparison} can be generalized under the proposed formalism. In all the following discussions, we assume a two-dimensional domain discretized by non-overlapping arbitrary polygonal volumes where the gradient of the scalar field needs to be reconstructed at the cell-centers. The extension of ideas presented herein to three dimensions is relatively straightforward.

\subsection {Green-Gauss schemes}
\label{GG}

The traditional Green-Gauss (GG) scheme can be obtained using the present formalism by defining the neighbourhood consisting of the face centers of the faces forming the cell ``$i$''. The geometric vectors are ${\bf a}_f = {\bf S}_f$ and ${\bf b}_f = {\bf x}_{f} - {\bf x}_i$. The area vector vector ${\bf S}_f = {\bf n}_f \Delta s_f$ where $\Delta s_f$ is the face area and ${\bf n}_f$ is the unit normal while ${\bf x}_f$ and ${\bf x}_i$ denote the position vectors of the face centers and the cell center respectively (See Fig.\ \ref{geo_cv}). The resulting dyadic-sum $\mathbb P$, for a two dimensional case, is given by,
\begin{equation}
\begin{split}
\mathbb P & =  \begin{bmatrix}
\displaystyle\sum_f n_x ({x}_f - {x}_i) \, \Delta s_f  & \displaystyle\sum_f n_x ({y}_f - {y}_i) \, \Delta s_f\\
\\
\displaystyle\sum_f n_y ({x}_f - {x}_i) \, \Delta s_f  & \displaystyle\sum_f n_y ({y}_f - {y}_i) \, \Delta s_f
\end{bmatrix} 
= \begin{bmatrix}
\Omega_i  & 0 \\
\\
0 & \Omega_i 
\end{bmatrix} ,
\end{split}
\label{eq:gg_dyadic_sum}
\end{equation}
where $(n_x,n_y)$ are the $x$ and $y$ components of $\mathbf{n}_f$ and $\Omega_i$ is the volume of the cell. The simplification of the terms in the matrix result from geometric identities (see Appendix \ref{secapp:gg_identities}). The dyadic-sum is therefore a diagonal matrix proportional to the identity tensor and can be inverted with ease. For this choice of vectors, the right hand side of Eq.\ \ref{eq5} simplifies to,
\begin{equation*}
\begin{split}
\sum_f {\bf a}_f~ \left({\bbf} \cdot {\nabla \phi_f}\right) 
&= \sum_f {\bf n}_f \Delta s_f \, \Big(({\bf x}_{f} - {\bf x}_i)\cdot {\nabla \phi_f}  \Big) = \sum_f {\bf n}_f \Delta s_f \, (\phi_{f} -\phi_i) .
\end{split}
\end{equation*}
Recognizing that $\sum_f \phi_i \, {\bf n}_f \Delta s_f = 0 $, the above expression substituted into the right hand side of Eq.\ \ref{eq5}, gives,
\begin{equation}
{\bf \nabla \phi}_i = \dfrac{1}{\Omega_i}\sum_f \phi_f  \, {\bf n}_f \Delta s_f ,
\label{ggeq}
\end{equation}
which is the Green-Gauss reconstruction formula. The face value $\phi_f$ is generally not available in cell-centered finite volume discretization and is typically obtained using suitable interpolation, which needs to be at least second order accurate to obtain consistent gradients \cite{syrakos2017critical}.


Besides the traditional Green-Gauss reconstruction, there exist other variants of the same family. One such scheme is the modified Green-Gauss (MGG) reconstruction which is derived from a different variant of the Gauss divergence theorem \cite{deka2018new}. This scheme, however, can be obtained through the proposed formalism by simply flipping the choices of the geometric vectors from GG reconstruction, i.e., ${\bf a}_f = {\bf x}_{f} - {\bf x}_i$ and ${\bf b}_f = {\bf S}_f = {\bf n}_f \Delta s_f $. The resulting dyadic sum is, therefore, just the transpose of $\mathbb{P}$ for GG, shown in Eq.\ \ref{eq:gg_dyadic_sum} and therefore enjoys the same ease of inversion. The directional derivative, however, is different from that of GG, leading to the right hand side of Eq.\ \ref{eq5} as,
\begin{equation}
\begin{split}
\sum_f {\bf a}_f~ \Big({\bf b} \cdot {\nabla \phi}\Big)_f 
&= \sum_f ({\bf x}_{f} - {\bf x}_i) \, \Big({\bf n}_f \Delta s_f  \cdot {\nabla \phi}  \Big) = \sum_f \frac{\partial \phi}{\partial n}\Big|_f \, ({\bf x}_{f} - {\bf x}_i) \, \Delta s_f .
\end{split}
\label{eq:mgg_rhs}
\end{equation}
Substituting this into the generalized formulation, we get the MGG reconstruction formula \cite{deka2018new},
\begin{equation}
{\bf \nabla \phi}_i = \dfrac{1}{\Omega}\sum_f \frac{\partial \phi}{\partial n}\Big|_f \, ({\bf x}_{f} - {\bf x}_i) \, \Delta s_f .
\label{mggeq}
\end{equation} 
For MGG, the face normal derivative can approximated on unstructured meshes using Zwart reconstruction \cite{deka2018new} or by central differences with a non-orthogonal correction \cite{jasak1996error}, which can essentially be interpreted as variants of the general idea of evaluating the directional derivative shown in Eq.\ \ref{eq7} in section \ref{sec2}. {It is instructive to note that, unlike traditional GG, the $\bbf$ vector for MGG does not pass through the cell center ``$i$'' except on genuinely orthogonal meshes. Therefore, an $O(h)$ accurate evaluation of the face normal derivative (which is necessary to obtain consistent gradients \cite{deka2018new}) typically requires the values and gradients of the scalar in the cell and its face-sharing neighbours. This makes MGG an implicit gradient reconstruction scheme.}


\subsection {Least squares schemes}
\label{lsq}

The least squares family of schemes are obtained by minimizing the error in the estimation of the values of a scalar in a specified stencil, in a least squares sense. The most commonly employed gradient scheme of the least squares family on unstructured meshes is the linear least squares reconstruction \cite{mavriplis2003revisiting,diskin2011comparison}, which boils down to minimizing the function $G^*$ defined by,
\begin{equation*}
G^* = \sum_{j=1}^{nb} \Big(\Delta \phi_j - \phi_{x,i} \Delta x_j - \phi_{y,i} \Delta y_j\Big)^2 ,
\end{equation*}
where, $\Delta x_j = x_j - x_i$, $\Delta y_j = y_j - y_i$, $\Delta \phi_j = \phi_j - \phi_i$, and $(\phi_{x,i},\phi_{y,i})$ are the $x$ and $y$ components of the $\gradphi_i$, respectively. The minimization problem reduces to a system of linear equations obtained by setting $\partial G^*/\partial \phi_{x,i} = 0$ and $\partial G^*/\partial \phi_{y,i} = 0$. Naturally due to the quadratic nature of the functional it is easy to see that the resulting system of equations contains a symmetric matrix that needs to be inverted to obtain the gradients. This property leads to the generalization in the proposed framework that least squares family of schemes correspond to a symmetric dyadic sum tensor, implying $\baf$ and $\bbf$ are co-linear. Therefore, the linear least squares scheme can be obtained from the generalized formalism, by first defining the neighbourhood at the centroids of the face-sharing cells\footnote{Extended stencils are also used which could include node-sharing cell centers or cell centers of cells sharing the face-sharing neighbours of ``$i$''.}, i.e. $f = j$ (see Fig.\ \ref{geo_cv}) and then defining, ${\bf a}_f =  {\bf b}_f = {\bf x}_j - {\bf x}_i$. The corresponding dyadic sum,
\begin{equation*}
\begin{split}
\mathbb P &=  \begin{bmatrix}
\displaystyle\sum_{j} \Delta x_j^2 & \displaystyle\sum_{j} \Delta x_j\Delta y_j  \\
\\
\displaystyle\sum_{j} \Delta y_j\Delta x_j & \displaystyle\sum_{j} \Delta y_j^2  \\
      \end{bmatrix},
\end{split}
\label{lsqeq2}
\end{equation*}
and the directional derivative,
\begin{equation*}
 \bbf \cdot \gradphi_f = (\xj - \xc)\cdot \gradphi_j = \phi_j - \phi_i  = \Delta \phi_j \:,
\end{equation*}
substituted into Eq.\ \ref{eq5} gives the linear least squares reconstruction formula,
\begin{equation}
\renewcommand\arraystretch{2}  
\begin{bmatrix}
\displaystyle\sum_{j} \Delta x_j^2  & \displaystyle\sum_{j} \Delta x_j \Delta y_j\\
\displaystyle\sum_{j} \Delta x_j \Delta y_j & \displaystyle\sum_{j} \Delta y_j^2 \\
\end{bmatrix} 
\begin{bmatrix}
\phi_{x,i} \\
\phi_{y,i} \\
\end{bmatrix} = 
\begin{bmatrix}
\displaystyle\sum_{j} \Delta \phi_j \Delta x_j\\
\displaystyle\sum_{j} \Delta \phi_j \Delta y_j\\
\end{bmatrix} .
\label{ulsq}
\end{equation}
The above system corresponds to the unweighted variant (ULSQ) of linear least squares reconstruction. The weighted linear least squares reconstruction (WLSQ) can be obtained from the proposed formalism by simply scaling the geometric vectors as ${\bf a}_j =  {\bf b}_j = w_j({\bf x}_j - {\bf x}_i)$, where the weights typically chosen such that it is inversely proportional to the cell-to-cell distance, $w_j = |{\bf x}_j - {\bf x}_i|^{-q}$ where $q \in [1,2]$. 

A different variant of the least squares family can be constructed by choosing the neighbourhood as the face centers of the cell ``$i$'' and thereby defining ${\bf a}_f = {\bf S}_f = {\bf n}_f \Delta s_f$ and ${\bf b}_f = {\bf n}_f$. This results in the symmetric dyadic-sum $\mathbb P$,
\begin{equation*}
\mathbb{P} = \begin{bmatrix}
\displaystyle\sum_f n_x^2 \Delta s_f & \displaystyle\sum_f n_x n_y \Delta s_f\\
\\
\displaystyle\sum_f n_y n_x \Delta s_f & \displaystyle\sum_f n_y^2 \Delta s_f\\
\end{bmatrix}.
\label{facelsq}
\end{equation*}
For this choice of geometric vectors we obtain,
\begin{equation*}
\begin{split}
\sum_f {\bf a}_f~ \Big({\bf b}_f \cdot {\nabla \phi_f}\Big) 
&= \sum_f {\bf S}_{f} \, \Big({\bf n}_f  \cdot {\nabla \phi}  \Big) = \sum_f \frac{\partial \phi}{\partial n}\Big|_f \, {\bf n}_{f} \, \Delta s_f .
\end{split}
\end{equation*}
The generalized formulation in Eq.\ \ref{eq5} then becomes, 
\begin{equation}
\renewcommand\arraystretch{2}  
\begin{bmatrix}
\displaystyle\sum_f n_x^2 \Delta s_f & \displaystyle\sum_f n_x n_y \Delta s_f \\
\displaystyle\sum_f n_y n_x \Delta s_f & \displaystyle\sum_f n_y^2 \Delta s_f \\
\end{bmatrix} 
\begin{bmatrix}
\phi_{x,i} \\
\phi_{y,i} \\
\end{bmatrix} = 
\begin{bmatrix}
\displaystyle\sum_f \dfrac{\partial \phi}{\partial n}\Big|_f n_x \Delta s_f \\
\displaystyle\sum_f \dfrac{\partial \phi}{\partial n}\Big|_f n_y \Delta s_f \\
\end{bmatrix} ,
\label{flsq2}
\end{equation}
where the summation is over all faces of the cell, $n_x$ and $n_y$ are the components of the unit normal to the face. This is the same system of linear algebraic equations obtained by minimizing the function
\begin{equation*}
G = \sum_f \Big(\dfrac{\partial \phi}{\partial n}\Big|_f - \dfrac{\partial\phi}{\partial x}\Big|_i n_x - \dfrac{\partial\phi}{\partial y}\Big|_i n_y\Big)^2 \Delta s_f,
\end{equation*}
which is also referred to as the face area weighted least squares reconstruction (FLSQ), first proposed in the context of energy-conserving schemes by Mahesh et al. \cite{mahesh2004numerical}. More recently, this method has been utilised to devise well-balanced solvers for incompressible multiphase flows \cite{ghods2013consistent,MANIK2018228} and also as a face-flux reconstruction technique \cite{weller2014non,AGUERRE2018135}. {Similar to MGG, the FLSQ scheme reconstructs gradients from the face normal derivative and is therefore an implicit scheme}.

\section{Hybrid schemes}

Besides the GG and LSQ schemes discussed in the previous section, there are some schemes which cannot be clearly categorised into either of these families but have essence of both. One such scheme is the recently proposed Taylor-Gauss (TG) gradients \cite{oxtoby2019family,SYRAKOS2023108}, which is obtained through an oblique projection and hence is not a minimization problem but leads to a least squares like system of equations. In our proposed framework, the same scheme can be obtained by first choosing the neighbourhood as the points of intersection of each face with the lines joining the cell center ``$i$'' to the cell center of the respective face-sharing neighbour cell. The geometric vectors are then chosen as, $\baf = \mathbf{S}_f$, and $\bbf = \xj - \xc$, where $\mathbf{S}_f$ is the face area vector. The resulting dyadic sum is,
\begin{equation*}
\mathbb P = \begin{bmatrix}
\displaystyle\sum_f n_x (x_j - x_i) \Delta s_f & \displaystyle\sum_f n_x (y_j - y_i) \Delta s_f\\
\\
\displaystyle\sum_f n_y (x_j - x_i) \Delta s_f & \displaystyle\sum_f n_y (y_j - y_i) \Delta s_f\\
\end{bmatrix}.
\label{dsg1}
\end{equation*}
With the directional derivative evaluated the same way as in the linear least squares scheme, we get, 
\begin{equation*}
\begin{split}
\sum_f {\bf a}_f~ \Big({\bf b}_f \cdot {\nabla \phi_f}\Big) = \sum_f {\bf S}_{f} \, \Big(({\bf x}_j - {\bf x}_i) \cdot {\nabla \phi}  \Big) = \sum_f (\phi_j - \phi_i) \, {\bf n}_{f} \, \Delta s_f ,
\end{split}
\label{dsg2}
\end{equation*}
and the corresponding system of equations to obtain the TG gradients as, 
\begin{equation}
\renewcommand\arraystretch{1.5}  
\begin{bmatrix}
\displaystyle\sum_f n_x (x_j - x_i) \Delta s_f & \displaystyle\sum_f n_x (y_j - y_i) \Delta s_f\\
\displaystyle\sum_f n_y (x_j - x_i) \Delta s_f & \displaystyle\sum_f n_y (y_j - y_i) \Delta s_f\\
\end{bmatrix} 
\begin{bmatrix}
\phi_{x,i} \\
\phi_{y,i} \\
\end{bmatrix} = 
\begin{bmatrix}
\displaystyle\sum_f (\phi_j - \phi_i) n_x \Delta s_f\\
\displaystyle\sum_f (\phi_j - \phi_i) n_y  \Delta s_f \\
\end{bmatrix} .
\label{dsg3}
\end{equation}
The weighted variants of Taylor-Gauss scheme \cite{SYRAKOS2023108} can be easily obtained by scaling the geometric vectors similar to the weighted variants of LSQ. 

It is remarked upon by the authors of \cite{oxtoby2019family} that the TG gradient schemes have a flavor of GG due to the use of face normals in reconstructing the centroidal gradients but with the additional ability to create weighted variants like an LSQ scheme. If we notice closely, to construct the TG scheme through the generalized formulation, we effectively borrowed $\baf$ and $\bbf$ from GG and LSQ schemes, respectively. Therefore, the scheme is, in some sense, a hybrid of GG and LSQ class of schemes, despite not being derivable from either routes clearly. This has been explained in \cite{SYRAKOS2023108} through their unification philosophy based on orthogonal and oblique projections. Within the context of our proposed formalism, the construction of TG gradients show a logical way of devising hybrid schemes out of traditional schemes, which is through different permutations of their respective geometric vectors. The efficacy of such a hybrid scheme can only be determined through numerical studies which is beyond the scope of this paper. Nevertheless, the proposed generalization provides an alternate route to creating nominally first order accurate hybrid gradient schemes where the advantages of the constituent schemes can possible be gained.

\section{Generalization of flexible gradient schemes}

Gradients of scalar fields are required during finite volume computation of fluid flow problems, where inconsistent gradients could lead to incorrect solutions. However, certain gradients reconstruction schemes like the linear least squares, despite being consistent, give rise to numerical instabilities on highly distorted meshes \cite{mavriplis2003revisiting,shima2013green}. Therefore, the efficacy of a gradient reconstruction scheme cannot be judged only in terms of its accuracy. For least squares reconstruction, the choice of the weights provide an indirect means to introduce robustness into the scheme. But a better idea is to directly introduce a damping parameter into the scheme, which was proposed in a flexible gradient scheme by Nishikawa \cite{nishikawa2021flexible}. The key idea in \cite{nishikawa2021flexible} is first to obtain gradients using conventional least squares approach and then use these gradients to compute the ``$\alpha$-damped'' gradients, which are then averaged to obtain the centroidal gradients. The parameter `$\alpha$' can be controlled to directly introduce damping without compromising the accuracy of the computed gradients. Additionally, it is also shown that for a specific value of `$\alpha$', the gradients become fourth order accurate on uniform Cartesian meshes. In this section, we propose a generalization of the flexible gradient idea using the generalized gradient reconstruction formalism of section 2. To motivate this, we first demonstrate the construction of a flexible variant of the Modified Green-Gauss (MGG) scheme \cite{deka2018new}. 
\subsection{A flexible variant of the MGG gradient scheme}
The Modified Green-Gauss (MGG) scheme of Deka et al. \cite{deka2018new} shown in Eq.\ \ref{mggeq} in section 3, reconstructs centroidal gradients from face normal derivatives, which are evaluated using Zwart reconstruction (see \cite{deka2018new}),
\begin{equation}
\frac{\partial \phi}{\partial n}\Big|_f \equiv \gradphi_f \cdot \nf = \alpha \left(\frac{\phi_j - \phi_i}{|\xj - \xc|}\right) + \left( \frac{\gradphi_i + \gradphi_j}{2}\right) \cdot (\nf - \alpha \rf) ,
\label{eq:flex_mgg_1}
\end{equation}
where $\nf$ and $\rf$ are unit vectors along the face normal and along the line joining cell center ``$i$'' to the face sharing neighbour cell center ``$j$'', respectively. In principle, Zwart reconstruction expresses the vector $\nf$ in $(\gradphi_f \cdot \nf)$ as, $\alpha \rf + (\nf - \alpha \rf)$. The `$\alpha \rf$' contribution can be written in terms of the centroidal values since $\rf$ is along the direction of the line joining the cell centers (first part on the right hand side of \ref{eq:flex_mgg_1}) whereas the `$(\nf - \alpha \rf)$' contribution is expressed in terms of an average of the centroidal gradients (second part on the right hand side of Eq.\ \ref{eq:flex_mgg_1}). While the right hand side of Eq.\ \ref{eq:flex_mgg_1} is directly employed to obtain MGG, it is possible to re-arrange it as, 
\begin{equation}
\frac{\partial \phi}{\partial n}\Big|_f = \underbrace{\left( \frac{\gradphi_i + \gradphi_j}{2}\right) \cdot \nf}_{\mbox{Consistent}} + \underbrace{\alpha \left(\frac{\phi_j - \phi_i}{|\xj - \xc|} - \left( \frac{\gradphi_i + \gradphi_j}{2}\right) \cdot \rf \right)}_{\mbox{Damping}}.
\label{eq:flex_mgg_2}
\end{equation}
Here, the first term on the right hand side could be interpreted as a consistent $O(h)$ accurate evaluation of the face normal derivative. The second term can be viewed as a damping component since it effectively represents a difference in the approximation of the quantity $(\gradphi_f \cdot \rf)$, weighted by the parameter $\alpha$. One could also re-write the second term upto $O(h)$ accuracy, as, 
$\alpha\left((\phi_R - \phi_L)/|\xj-\xc|\right)$,
where, $\phi_L = \phi_i + \gradphi_i \cdot (\xf - \xc)$ and $\phi_R = \phi_j + \gradphi_j \cdot (\xf - \xj)$, which is equivalent in form to the damping component of the flexible scheme proposed by Nishikawa \cite{nishikawa2021flexible}. 
In the original formulation of MGG \cite{deka2018new}, the parameter $\alpha$ is a fixed constant, chosen as $(\nf\cdot\rf)$. This is a geometrically intuitive choice as it makes the vector `$(\nf - \alpha \rf)$' orthogonal to `$\nf$', thereby decreasing the implicit contribution with decreasing non-orthogonality. 
However, from a reconstruction accuracy perspective, the value of $\alpha$ is not constrained by any geometric considerations. Therefore, adopting $\alpha$ as a constant free-parameter, one can turn MGG into a flexible implicit gradient scheme. 
 
Instead of the implicit formulation, we could also re-frame the scheme into a two-step explicit procedure (similar to \cite{nishikawa2021flexible}) with the purpose of numerical stabilization. 
The weighted LSQ gradients are known to be consistent but unstable for Euler computations on distorted meshes \cite{mavriplis2003revisiting,shima2013green}. We therefore compute these gradients as the first step which can be then employed in determining the normal derivative in Eq.\ \ref{eq:flex_mgg_1} in the second step. This makes the process explicit and use of the damping parameter provides the flexibility to ensure convergence of the numerical approach with the gradients obtained by this two-step approach. While we have introduced this concept with reference to the MGG gradients, it is indeed possible to generalize this approach as discussed in the following section.

\subsection{A general flexible scheme}
For a general gradient scheme, it is shown through the proposed formalism in section 2 that, the centroidal gradients are reconstructed from the quantity $(\bbf \cdot \gradphi_f)$ in a defined neighbourhood. Similar to the face normal splitting in Zwart reconstruction, we can split the vector $\bbf$ as,
\begin{equation*}
 \frac{\bbf}{|\bbf|} = \alpha \rf + \left(\frac{\bbf}{|\bbf|} - \alpha \rf\right),
\end{equation*}
which gives,
\begin{equation}
\begin{aligned}
   {\bf b}_f \cdot \nabla {\phi}_f &= \alpha |\bbf|\left(\rf \cdot \nabla {\phi}_f\right) + ({\bf b}_f - \alpha |\bbf| \rf)\cdot \gradphi_f \\
	&= \alpha |\bbf| \left(\frac{\phi_j - \phi_i}{|\xj - \xc|}\right) + (\bbf - \alpha|\bbf|\rf ) \cdot \left( \frac{{\nabla \phi}_j + {\nabla \phi}_i}{2} \right).
   \end{aligned}
   \label{eq:alpha_eq16}
\end{equation}
This can again be shown as a sum of a consistent and damping part similar to Eq.\ \ref{eq:flex_mgg_2}. Therefore, by substituting Eq.\ \ref{eq:alpha_eq16} into the right hand side of Eq.\ \ref{eq5}, we can obtain a nominally first order accurate implicit flexible variant of any gradient scheme\footnote{From error analyses similar to the ones outlined in appendix \ref{secapp:error_terms}, it is easy to show that the generalized $\alpha$-damped gradients in Eq.\ \ref{eq:alpha_eq16} reconstruct nominally first order accurate centroidal gradients.}. 

The generalized flexible scheme can also be adopted into a two-step explicit stabilization mechanism for an existing `unstable' scheme similar the one proposed by Nishikawa \cite{nishikawa2021flexible}. Suppose we have a gradient scheme represented by the geometric vectors $(\baf,\bbf)$ that is unstable on a poorly conditioned mesh. The idea would be to compute the centroidal gradients using that scheme as the first step (lets denote them as $\tilde{\gradphi}$). In the second step, we substitute these gradients into the implicit part on the right hand side of Eq.\ \ref{eq:alpha_eq16}, i.e,
\begin{equation}
\begin{aligned}
   {\bf b}_f \cdot \nabla {\phi}_f = \alpha |\bbf| \left(\frac{\phi_j - \phi_i}{|\xj - \xc|}\right) + (\bbf - \alpha|\bbf|\rf ) \cdot \left( \frac{\tilde{{\nabla \phi}_j} + \tilde{{\nabla \phi}_i}}{2} \right).
   \end{aligned}
   \label{eq:alpha_eq17}
\end{equation}
This substituted into Eq.\ \ref{eq5} leads to the two-step flexible variant of a general scheme.
For numerical stabilization, the value of $\alpha$ has to be suitably chosen.
It is, however, not possible for the generalization to specify how the flexible parameter has to be chosen to increase the robustness of a scheme, especially on arbitrary unstructured meshes. One could perform a modified wave-number analysis on a one-dimensional uniform grid to get an idea of the role of $\alpha$, similar to Nishikawa \cite{nishikawa2021flexible}. But it is easy to show, that for a one-dimensional grid the geometric vectors degenerate to constant scalars and therefore the centroidal gradients degenerate to the same expression for all flexible schemes, including the one shown in \cite{nishikawa2021flexible}. Therefore, one can expect that on meshes closer to a uniform Cartesian mesh, decreasing $\alpha$ from unity should increase the numerical damping (at $\alpha = 1$ the undamped scheme is recovered).
We, however, stress that the trends need not necessarily reciprocate in the same manner on arbitrary unstructured meshes and can only be determined by numerical tests. Nonetheless, the presence of a flexible parameter does provide a means to introduce damping for highly distorted meshes where gradients computed from a consistent scheme could cause numerical instabilities in flow computations.
 

\section{General discussions and summary}
\label{sec5}

In this work, we have presented a generalization of gradient reconstruction schemes in cell-centered unstructured finite volume method. The key components of the proposed framework are a pair of geometric vectors defined in a suitably chosen neighbourhood of a cell, different choices of which yield different gradient schemes. Commonly employed gradient reconstruction schemes generally belong to either the Green-Gauss family or the least squares family, both of which are built from starkly different philosophies. The proposed formalism shows that they can be mathematically unified under a common framework. But as these schemes can be conveniently derived from their own mathematical principles, it might appear that the proposed generalization merely serves as a post facto mathematical `fitting' tool. The authors, however, believe that it is not just the intended purpose of this framework. 

The generalization provides some important insights into the basic philosophy of gradient reconstruction. It shows that every nominally first order accurate gradient scheme effectively reconstructs centroidal gradients from the directional derivative of the scalar along specific directions in a defined neighbourhood. This viewpoint is not trivial and, to the best of the authors' knowledge, has not been perceived before, as it is not seen from either the divergence theorem or the minimization based route. The only known generalization of gradient schemes is based on the idea of orthogonal and oblique projections proposed by Syrakos et al.\ \cite{SYRAKOS2023108}. Their framework has a similar flavor to the one proposed here as it also has a dyadic sum that requires inversion. But there is an important distinction between the two formalisms that is worth understanding. The framework proposed by Syrakos et al.\ \cite{SYRAKOS2023108} begins with obtaining a first order accurate estimate of the change in the scalar along a particular direction, using the value and gradients at the centroid of a cell ``$i$'', i.e., 
\begin{equation}
 \Delta \phi_f = \phi_f - \phi_i = \mathbf{R}_f \cdot \gradphi_i  .
\label{eq:syrakos_1}
\end{equation}
Here, ``$f$'' denotes a point in a compact neighbourhood about ``$i$''. The above expression is now multiplied by another vector $\mathbf{V}_f$, and then summed across the neighbourhood to obtain the generalized gradient formula,
\begin{equation}
\sum_f \mathbf{V}_f \Delta \phi_f = \Big(\sum_f \mathbf{V}_f \otimes \mathbf{R}_f \Big) \gradphi_i  .
\label{eq:syrakos_2}
\end{equation}
Comparing this to Eq.\ \ref{eq5}, it can be easily seen that the vector $\mathbf{V}_f$ is equivalent to $\baf$. Also, $\Delta \phi_f$ on the left hand side of Eq.\ \ref{eq:syrakos_2} can be seen as an evaluation of the quantity $\mathbf{R}_f \cdot \gradphi_i$ (from Eq.\ \ref{eq:syrakos_1}), which is effectively the directional derivative of $\phi$ along the vector $\mathbf{R}_f$. Therefore, we can draw an equivalence of $\mathbf{R}_f$ to $\bbf$, but that is not strictly true. The way in which the formulation proposed by Syrakos et al.\ \cite{SYRAKOS2023108} is built, $\mathbf{R}_f$ always points from the cell center ``$i$'' to some point ``$f$'' in the defined neighbourhood. The location of ``$f$'' is kept arbitrary, which provides a degree of freedom in the choice of $\mathbf{R}_f$, but the vector always passes through the cell center ``$i$''. This restriction is not imposed on the vector $\bbf$ in our framework. Physically it means that at conception, the framework proposed by Syrakos et al.\ \cite{SYRAKOS2023108} reconstructs centroidal gradients from the variation of the scalar between the cell center and an arbitrary point in the neighbourhood, whereas we propose that the same variation could be taken between any two arbitrary points in a compact neighbourhood. Of course, it is always possible to express an arbitrary $\bbf$ vector as a sum of $\mathbf{R}_f$ and another vector to obtain some form of a generalization like \cite{SYRAKOS2023108}, but at a philosophical level, our proposition has a different interpretation. It says that instead of the $\phi$ difference between the cell center and an arbitrary `point', it is, in fact, the change in $\phi$ along an arbitrary `direction' that is required to reconstruct nominally first order accurate centroidal gradients. 

The implication of this insight is apparent from the re-derivation of some of the existing gradient schemes. From section 3, one can notice that for the traditional Green-Gauss and least squares reconstruction, the vector $\bbf$ is either $\xf - \xc$ (for GG) or $\xj - \xc$ (for LSQ), both of which pass through ``$i$'', and therefore can be comfortably derived using the framework of Syrakos et al.\ \cite{SYRAKOS2023108} as well as ours. But for schemes like the MGG or FLSQ, where $\bbf$ is along the face normal (which for an arbitrary unstructured mesh does not pass through cell centers), it not clear how the framework of Syrakos et al.\ \cite{SYRAKOS2023108} can be used to derive them. The arbitrariness of $\bbf$ also allows for construction of hybrid schemes by just using $\baf$ and $\bbf$ vectors directly from any two known schemes as shown in section 4. In addition, the perspective of reconstructing centroidal gradients from directional derivatives allows us to create flexible variants of known schemes as shown in section 5. The idea is to introduce an additional damping component in the evaluation of the directional derivative that is dependent on the centroidal gradients. These centroidal gradients could be evaluated apriori from a base scheme that is possibly unstable and then substituted into the $\alpha$-damped version of the directional derivatives. With a suitably chosen value of the parameter $\alpha$, the reconstructed gradients can introduce stability to finite volume computations without compromising the gradient accuracy. From a user perspective, this provides a simple route to increasing the robustness of their solver by modifying their choice of gradient reconstruction scheme into the flexible variant.

While the novelties of the proposed framework have been explained, it is incumbent upon us to highlight some of its shortcomings. The mathematical generalization proposed here does not outline a clear directive as to how the geometric vectors are to be chosen to construct a new gradient reconstruction scheme. The only constraints on the geometric vectors are,
\begin{enumerate}
\item The geometric vectors should depend on the mesh metrics, 
\item Their magnitudes should be $O(h^k)$, where $k>0$, and,
\item The sum of their dyadic product over the chosen neighbourhood should be an invertible tensor. 
\end{enumerate}
These guidelines are reasonably arbitrary from the point of view of directing possible choices of $\baf$ and $\bbf$ to construct a scheme. The generalization, however, does indicate that one need not necessarily follow the traditional routes, like the divergence theorem or a minimization idea, to derive a first order accurate gradient reconstruction scheme. A good example is the Taylor-Gauss (TG) gradient scheme \cite{oxtoby2019family,SYRAKOS2023108} which is derivable neither from the GG or the LSQ route, but can be easily obtained from the unified framework of Syrakos et al.\ \cite{SYRAKOS2023108} as well as ours. Despite not being derivable from either routes, the scheme has flavor of both GG and LSQ families \cite{oxtoby2019family,SYRAKOS2023108}, simply because the geometric vectors that define the scheme are borrowed from GG or LSQ. The generalization, therefore, allows us to see the mathematical hybridization that exists in the TG scheme, and could possibly pave a way for more such hybrid schemes in the future. Nevertheless, the framework does not specify how such a hybridization could be `better' or `worse' that a conventional scheme. As like in most cases, these can only be determined by numerical tests. Similarly in the extension to flexible gradient schemes, the generalization does not specify how to choose the flexible parameter to increase robustness of a given scheme. It is not generally possible to do so on arbitrary unstructured meshes as the degree of damping is highly dependent on the scheme as well as the solution. This is demonstrated in a future work by the same authors, where the proposed generalization of $\alpha$-damped scheme is used to construct flexible variants of some known schemes.

\backmatter

\bmhead{Acknowledgments}

The authors would like to express their gratitude to Thibault Pringuey for his write-up on CFD Online detailing the ideas behind {\texttt {fvc::reconstruct}} in OpenFOAM which motivated the present study. The last author would like to acknowledge the support from Science and Engineering Research Board (SERB), Government of India through the MATRICS grant MTR/2019/000241 during the course of this work.

\section*{Statements and Declarations}

\subsection*{Funding}
This work was supported by Science and Engineering Research Board (SERB), Government of India through the MATRICS
grant MTR/2019/000241.

\subsection*{Competing interests}
The authors declare that they have no competing interests.

\subsection*{Author contributions}
All authors contributed equally to the work.

\subsection*{Availability of data and materials}
No data were generated or analyzed in this study.


\begin{appendices}

\section{Error analysis}
\label{secapp:error_terms}

\subsection{Term 1:}
\label{secapp:error_term_1}

In obtaining Eq.\ \ref{eq5} from Eq.\ \ref{eq3}, the following term is ignored, 
\begin{equation}
 E_1 = -\mathbb{P}^{-1} \sum_f (\baf \otimes \bbf) \cdot ({\nabla \gradphi}_c \cdot (\xf - \xc) + \epsilon_1).
\end{equation}
Here, $\epsilon_1 \sim O(h^2)$. From Cauchy-Schwarz and triangle inequalities, we can write,
\begin{equation*}
        ||E_1|| \leq ||{\mathbb{P}}^{-1}|| \sum_f \left(||\baf \otimes \bbf||\right) \left(||(\xf - \xc)||\right) \left(||{\nabla \gradphi}_c||\right),
\end{equation*}
where, norms are defined in the usual sense of dot product for vectors and tensor inner product for second order tensors. From equation \ref{eq0}, we can write $||\baf \otimes \bbf|| \sim {O}(h^{p+q})$ and therefore, $||{\mathbb{P}}^{-1}|| \sim {O}(h^{-p-q)})$. Assuming a reasonably compact neighbour stencil, $||(\xf - \xc)|| \sim {O}(h)$. Since any norm of the scalar field $\phi$ is independent of the grid metrics, we can therefore write, 
    \begin{equation*}
        ||E_1|| \sim {O}(h).
    \end{equation*}

\subsection{Term 2:}
\label{secapp:error_term_2}

The second term on the right hand side of Eq.\ \ref{eq7} is an error term, which when substituted into Eq.\ \ref{eq5} gives the contribution,
\begin{equation}
E_2 = \mathbb{P}^{-1} \sum_f \baf \left( (\delftwo + \delfone) ( {\nabla \nabla \phi}_f : (\bbf \otimes \bbf ) + \epsilon_2) \right).
\end{equation}
Here, $\epsilon_2 \sim O(h^2)$.
Similar to $E_1$, the norm of $E_2$ can be written as,
    \begin{equation*}
    \begin{aligned}
        ||E_2|| & \leq ||\mathbb P^{-1} || \sum_f ||{\bf a}_f|| \: ( |\delftwo + |\delfone|) \: (||{{\bf b}_f} \otimes {{\bf b}_f}|| ) .
    \end{aligned}
    \end{equation*}
    The choice of $\delta$ were such that $|\delfone \bbf| \sim O(h)$ and $|\delftwo \bbf| \sim O(h)$. Therefore, from Eq.\ \ref{eq0}, we can write, $|\delfone| \sim O(h^{1-q})$ and  $|\delftwo| \sim O(h^{1-q})$. With $||\bbf \otimes \bbf|| \sim O(h^{2q})$, and $||\baf|| \sim O(h^p)$, while $||\mathbb{P}^{-1}|| \sim {O}(h^{-p-q})$, we get, 
     \begin{equation*}
        ||E_2|| \sim {O}(h).
    \end{equation*}

\subsection{Term 3:}
\label{secapp:error_term_3}

The values of $\phi$ at $\xfone$ and $\xftwo$ can be interpolated from neighbour cells. Let us consider an interpolation that is $O(h^r)$ accurate. In that case, the error term incurred in Eq.\ \ref{eq9} due to interpolation can be written as,
\begin{equation}
 E_3 = \mathbb{P}^{-1} \sum_f \baf \left( \dfrac{\epsilon_{\phi_f,2} - \epsilon_{\phi_f,1}}{\delftwo - \delfone} \right),
\end{equation}
where $\epsilon_{\phi_f,1}, \epsilon_{\phi_f,2} \sim O(h^{r})$, are the interpolation errors for $\phi_{f,1}, \phi_{f,2}$. Using Cauchy-Schwartz and triangle inequalities, we can write,
\begin{equation*}
    \begin{aligned}
       ||E_3|| & \leq  ||\mathbb P^{-1}|| \sum_f ||{\bf a}_f|| \: (|| \epsilon_{\phi_f,2} - \epsilon_{\phi_f,1}||) \: \Big| \Big|\dfrac{1}{{(\delta^{(2)}_f - \delta^{(1)}_f)}}\Big| \Big| .
    \end{aligned}
    \end{equation*}
    Since $|\delfone \bbf| \sim O(h)$ and $|\delftwo \bbf| \sim O(h)$, $||\delta^{(2)}_f - \delta^{(1)}_f||^{-1}\sim O(h^{q-1})$. Since, $||\baf|| \sim O(h^p)$, while $||\mathbb{P}^{-1}|| \sim {O}(h^{-p-q})$, we get 
    \begin{equation*}
        ||E_3|| \approx {\bf O}(h^{r-1})
    \end{equation*}
    Therefore, for $E_3$ to be at max $O(h)$, it is required that $r \geq 2$.

\section{Some useful identities}
 \label{secapp:gg_identities}

The dyadic sum defined for the Green-Gauss family (GG and MGG reconstructions) can be simplified employing a simple identity for the volume of a cell. While this identity is commonplace in unstructured computations, we present its concise derivation for the sake of completeness and benefit of interested readers. The Gauss divergence theorem for a vector field $\bf F$ reads,
 \begin{equation}
 \int_\Omega \nabla \cdot {\bf F} \, \ud\Omega  = \oint {\bf F} \cdot {\bf n} \, \ud s .
 \end{equation}
Letting $\bf F = {\bf c}\phi$ where $\phi$ is a scalar and $\bf c$ is a constant vector gives (after simplification),
 \begin{equation}
 \int_\Omega \nabla \phi \, \ud\Omega  = \oint \phi \, {\bf n} \, \ud s .
 \label{gdteq}
 \end{equation}
\noindent Consider that $\phi = \mathbf{x} \cdot \mathbf{e}$, where $\mathbf{x}$ is a position vector and $\mathbf{e}$ is any unit vector in Cartesian co-ordinates. Then, Eq.~\ref{gdteq} becomes,
   \begin{equation}
  \int_\Omega \nabla (\mathbf{x}\cdot \mathbf{e}) \, \ud\Omega = \oint (\mathbf{x}\cdot \mathbf{e}) \, {\bf n}\ud s = \sum_f (\mathbf{x}_f \cdot \mathbf{e}_f)\, \mathbf{n}_f \Delta s_f .
  \label{gdteq2}
 \end{equation} 
This is an exact relation which computes the integral over each face using single-point Gauss quadrature. One can then show with little effort that Eq. \ref{gdteq2} reduces to the following relations if ${\bf e} = [1~0]^T$,
\begin{eqnarray*}
&\displaystyle\sum_f(x_f-x_i) \, n_x \, \Delta s_f =  \Omega_i ,& \\
&\displaystyle\sum_f(x_f-x_i) \, n_y \, \Delta s_f =  0 .&
\end{eqnarray*} 
A similar exercise for ${\bf e} = [0~1]^T$ gives,
\begin{eqnarray*}
&\displaystyle\sum_f(y_f-y_i) \, n_x  \, \Delta s_f =  0 ,& \\
&\displaystyle\sum_f (y_f-y_i) \, n_y  \, \Delta s_f =  \Omega_i ,&
\end{eqnarray*}
where $\Omega_i$ is the volume of the cell $i$, the summation is over the faces of $i$ and the identity $\sum_f {\bf n}_f \Delta s_f  = 0$ has been employed. The latter identity itself can be easily obtained from Eq. \ref{gdteq} by setting $\phi$ equal to any constant scalar. The three dimensional analogues of all the above identities follow easily.

\end{appendices}

\bibliography{sample.bib}   

\end{document}